\documentclass[11pt]{article}
\usepackage{diagrams}
\newcommand{\dated}{\mbox{} \hfill {\small [{\tt \today}]}} 
\newtheorem{theorem}{Theorem}[section]
\newtheorem{lemma}[theorem]{Lemma}
\newtheorem{corollary}[theorem]{Corollary}
\newtheorem{proposition}[theorem]{Proposition}
\newtheorem{df}[theorem]{Definition}
\newenvironment{definition}{\begin{df} \rm}{\end{df}}

 \usepackage{amsmath,amssymb,amscd}
\newcommand{\pf}[1]{\trivlist \item[\hskip\labelsep\it #1\ ]}
\newcommand{\varpf}[1]{\trivlist \item[\hskip\labelsep\sc #1:]}
\newcommand{\qedbox}{$\rlap{$\sqcap$}\sqcup$}
\newcommand{\qed}{\qquad \qedbox \endtrivlist}
\newcommand{\varqed}{\hfill \rule{0.6em}{0.6em} \endtrivlist}
\newenvironment{proof}{\pf{Proof}}{\qed}

\newenvironment{remark}{\pf{Remark}}{\endtrivlist}
\newenvironment{remarks}{\pf{Remarks} 
   \begin{enumerate}}{\end{enumerate} \endtrivlist}
\newenvironment{example}{\pf{Example}}{\endtrivlist}
\newenvironment{examples}{\pf{Examples} 
   \begin{enumerate}}{\end{enumerate} \endtrivlist}
\newenvironment{items}{
  \begin{enumerate} 
                    
  }{\end{enumerate}}
\newenvironment{alphitems}{
  \begin{enumerate} 
                    
  }{\end{enumerate}}
\newenvironment{keywords}{\noindent\small {\it Keywords\/}:}{\vskip 4pt}
\newenvironment{classification}{\noindent\small 2000 {\it Mathematics Subject
Classification\/}:}{\vskip 12pt}

\newcommand{\comps}{{\mathbb C}}

\newcommand{\ints}{{\mathbb Z}}
\newcommand{\posints}{{\mathbb N}}

\newcommand{\tensor}{\otimes}

\newcommand{\Tensor}{\hat{\otimes}}

\newcommand{\cstar}{{C^\ast}}

\newcommand{\A}{{\mathfrak A}}
\newcommand{\B}{{\mathfrak B}}

\newcommand{\AM}{\operatorname{AM}}
\newcommand{\WAP}{\operatorname{WAP}}
\newcommand{\wWC}{\operatorname{{\sigma}WC}}
\title{Dual Banach algebras: Connes-amenability, \\
normal, virtual diagonals, and \\
injectivity of the predual bimodule}
\author{{\it Volker Runde}\thanks{Research supported by NSERC under grant no.\ 227043-00.}}
\date{}
\begin{document}
\maketitle
\begin{abstract}
Let $\A$ be a dual Banach algebra with predual $\A_\ast$ and consider the following assertions: (A) $\A$ is Connes-amenable; (B) $\A$ has a normal, virtual diagonal; (C) $\A_\ast$ is an injective $\A$-bimodule.
For general $\A$, all that is known is that (B) implies (A) whereas, for von Neumann algebras, (A), (B), and (C) are equivalent. We show that (C) always implies (B) whereas the converse is false for $\A = M(G)$ where $G$ is
an infinite, locally compact group. Furthermore, we present partial solutions towards a characterization of (A) and (B) for $\A = B(G)$ in terms of $G$: For amenable, discrete $G$ as well as for certain compact $G$, they are
equivalent to $G$ having an abelian subgroup of finite index. The question of whether or not (A) and (B) are always equivalent remains open. However, we introduce a modified definition of a normal, virtual diagonal and, using 
this modified definition, characterize the Connes-amenable, dual Banach algebras through the existence of an appropriate notion of virtual diagonal.
\end{abstract}
\section*{Introduction}
In \cite{JKR}, B.\ E.\ Johnson, R.\ V.\ Kadison, and J.\ Ringrose introduced a notion of amenability for von Neumann algebras which modifies Johnson's original definition for general Banach algebras (\cite{Joh1}) in the
sense that it takes the dual space structure of a von Neumann algebra into account. This notion of amenability was later dubbed Connes-amenability by A.\ Ya.\ Helemski\u{\i} (\cite{HelSb}).
\par
In \cite{Run}, the author extended the notion of Connes-amenability to the larger class of dual Banach algebras (a Banach algebra is called dual if it is a dual Banach space such that multiplication is separately
$w^\ast$-continuous). Examples of dual Banach algebras (besides von Neumann algebras) are, for example, the measure algebras $M(G)$ of locally compact groups $G$. In \cite{Run2}, the author proved that a locally compact 
group $G$ is amenable if and only if $M(G)$ is Connes-amenable --- thus showing that the notion of Connes-amenability is of interest also outside the framework of von Neumann algebras.
\par
In \cite{Eff}, E.\ G.\ Effros showed that a von Neumann algebra is Connes-amenable if and only if it has a so-called normal, virtual diagonal. Like Connes-amenability, the notion of a normal, virtual diagonal
adapts naturally to the context of general dual Banach algebras. It is not hard to see that a dual Banach algebra with a normal, virtual diagonal is Connes-amenable (the argument from the von Neumann algebra case
carries over almost verbatim; see \cite{CG}).
\par
Let $\A$ be a dual Banach algebra with (not necessarily unique) predual $\A_\ast$; it is easy to see that $\A_\ast$ is a closed submodule of $\A^\ast$. Consider the following three statements:
\begin{itemize}
\item[(A)] $\A$ is Connes-amenable.
\item[(B)] $\A$ has a normal, virtual diagonal.
\item[(C)] $\A_\ast$ is an injective $\A$-bimodule in the sense of \cite{Hel}.
\end{itemize}
\par
If $\A$ is a von Neumann algebra, then (A), (B), and (C) are equivalent (the equivalence of (A) and (B) was mentioned before; that they are equivalent to (C) is proved in \cite{HelSb}). If $\A = M(G)$ for a locally compact 
group $G$, then (A) and (B) are also equivalent (\cite{Run3}). For a general dual Banach algebra $\A$, we know that (B) implies (A), but nothing else seems to be known about the relations between (A), (B), and (C).
\par
As we shall see in the present paper, (C) always implies (B) --- and thus (A) --- whereas the converse need not hold in general: this answers a question by A.\ Ya.\ Helemski\u{\i}
(\cite[Problem 24]{LoA}) in the negative. The counterexample is the measure algebra $M(G)$ for any infinite, amenable, locally compact group $G$; the proof relies on recent work by H.\ G.\ Dales and M.\ Polyakov
(\cite{DP}). (As O.\ Yu.\ Aristov informed us upon seeing a preprint version of this paper, it had previously been shown by S.\ Tabaldyev that the Banach $\ell^1(\ints)$-bimodule $c_0(\ints)$ 
is not injective, which already answers Helemski\u{\i}'s question; see \cite{Tab}.)
\par
The Fourier--Stieltjes algebra $B(G)$ of a locally compact group $G$, as introduced in \cite{Eym}, is another example of a dual Banach algebra. In view of \cite{Run2}, \cite{Run3}, and \cite{Run4}, it is not farfetched
to conjecture that (A) and (B) for $\A = B(G)$ are equivalent and hold true if and only if $G$ has an abelian subgroup of finite index. Even though we are not able to settle this conjecture in full generality, we can corroborate it for 
certain $G$: (A) and (B) hold for $B(G)$ --- with $G$ discrete and amenable or a topological product of finite groups --- if and only if $G$ has an abelian subgroup of finite index.
\par
In the last section of the paper we modify the definition of a normal, virtual diagonal by introducing what we call a $\wWC$-virtual diagonal. For a dual Banach algebra $\A$, we then consider the statement:
\begin{itemize}
\item[(B')] $\A$ has a $\wWC$-virtual diagonal.
\end{itemize}
Unlike for (A) and (B), we can show that (A) and (B') are indeed equivalent. It thus seems that the notion of a $\wWC$-virtual diagonal seems to be the more natural one to consider in the context of 
Connes-amenability if compared with the notion of a normal, virtual diagonal.
\section{Preliminaries}
\subsection{Notions of amenability}
We start with the definition of a dual Banach module:
\begin{definition} \label{dualdef}
Let $\A$ be a Banach algebra. A Banach $\A$-bimodule $E$ is called {\it dual\/} if it is the dual of some Banach space $E_\ast$ such that, for each $a \in \A$, the maps
\[
  E \to E, \quad x \mapsto \left\{ \begin{array}{c} a \cdot x, \\ x \cdot a \end{array} \right.
\]
are $\sigma(E,E_\ast)$-continuous.
\end{definition}
\begin{remarks}
\item The predual space $E_\ast$ in Definition \ref{dualdef} need not be unique. Nevertheless, $E_\ast$ will always be clear from the context, so that we can speak of the $w^\ast$-topology on $E$ without ambiguity.
\item It is easily seen that a dual Banach space $E$ (with predual $E_\ast$) which is also a Banach $\A$-bimodule is a dual Banach $\A$-bimodule if and only if $E_\ast$ is a closed submodule of $E^\ast$. Hence,
our definition of a dual Banach $\A$-bimodule coincides with the usual one (given in \cite{LoA}, for instance).
\end{remarks}
\par
Let $\A$ be a Banach algebra, and let $E$ be a Banach $\A$-bimodule. A {\it derivation\/} from $\A$ to $E$ is a bounded, linear map $D \!: \A \to E$ satisfying
\[
  D(ab) = a \cdot Db + (Da) \cdot b \qquad (a,b \in \A).
\]
A derivation $D \!: \A \to E$ is called {\it inner\/} if there is $x \in E$ such that
\[
  Da = a \cdot x - x \cdot a \qquad (x \in \A).
\]
\begin{definition} \label{amdef}
A Banach algebra $\A$ is called {\it amenable\/} if every derivation from $\A$ into a dual Banach $\A$-bimodule is inner.
\end{definition}
\par
The terminology is, of course, motivated by \cite[Theorem 2.5]{Joh1}: A locally compact group $G$ is amenable if and only if its group algebra $L^1(G)$ is amenable in the sense of Definition \ref{amdef}.
\par
For some classes of Banach algebra, Definition \ref{amdef} seems to be ``too strong'' in the sense that it only characterizes fairly uninteresting examples in those classes: A von Neumann algebra is amenable if and only if it is
subhomogeneous (\cite{Was}), and the measure algebra $M(G)$ of a locally compact group $G$ is amenable if and only if $G$ is discrete and amenable (\cite{DGH}). 
\par
Both von Neumann algebras and measure algebras are dual Banach algebras in the sense of the following definition:
\begin{definition}
A Banach algebra $\A$ which is a dual Banach $\A$-bimodule is called a {\it dual Banach algebra\/}.
\end{definition}
\begin{examples}
\item Every von Neumann algebra is a dual Banach algebra.
\item The measure algebra $M(G)$ of a locally compact group $G$ is a dual Banach algebra (with predual ${\cal C}_0(G))$.
\item If $E$ is a reflexive Banach space, then ${\cal B}(E)$ is a dual Banach algebra (with predual $E \Tensor E^\ast$, where $\Tensor$ denotes the projective tensor product of Banach spaces).
\item The bidual of every Arens regular Banach algebra is a dual Banach algebra.
\end{examples}
\par
We shall now introduce a variant of Definition \ref{amdef} for dual Banach algebras that takes the dual space structure into account:
\begin{definition}
Let $\A$ be a dual Banach algebra, and let $E$ be a dual Banach $\A$-bimodule. An element $x \in E$ is called {\it normal\/} if the maps
\[
  \A \to E, \quad a \mapsto \left\{ \begin{array}{c} a \cdot x, \\ x \cdot a \end{array} \right.
\]
are $w^\ast$-continuous. The set of all normal elements in $E$ is denoted by $E_\sigma$. We say that $E$ is {\it normal\/} if $E = E_\sigma$.
\end{definition}
\begin{remark}
It is easy to see that, for any dual Banach $\A$-bimodule $E$, the set $E_\sigma$ is a norm closed submodule of $E$. Generally, however, there is no need for $E_\sigma$ to be $w^\ast$-closed.
\end{remark}
\begin{definition} \label{condef}
A dual Banach algebra $\A$ is called {\it Connes-amenable\/} if every $w^\ast$-con\-ti\-nu\-ous derivation from $\A$ into a normal, dual Banach $\A$-bimodule is inner.
\end{definition}
\begin{remarks}
\item ``Connes''-amenability was introduced by B.\ E.\ Johnson, R.\ V.\ Kadison, and J.\ Ringrose for von Neumann algebras in \cite{JKR}. The name ``Connes-amenability'' seems to originate in \cite{HelSb}, probably in reverence towards 
A.\ Connes' fundamental paper \cite{Con1}.
\item For a von Neumann algebra, Connes-amenability is equivalent to a number of important properties, such as injectivity and semidiscreteness; see \cite[Chapter 6]{LoA} for a relatively self-contained account.
\item The measure algebra $M(G)$ of a locally compact group $G$ is Connes-amenable if and only if $G$ is amenable (\cite{Run2}).
\end{remarks}
\subsection{Virtual diagonals}
Let $\A$ be a Banach algebra. Then $\A \Tensor \A$ is a Banach $\A$-bimodule via
\[
  a \cdot (x \tensor y) := ax \tensor y \quad\text{and}\quad (x \tensor y) \cdot a := x \tensor ya \qquad (a,x,y \in \A),
\]
so that the multiplication map
\[
  \Delta \!: \A \Tensor \A \to \A, \quad a \tensor b \mapsto ab
\]
becomes a homomorphism of Banach $\A$-bimodules.
\par
The following definition is also due to B.\ E.\ Johnson (\cite{Joh2}):
\begin{definition} \label{diagdef}
A {\it virtual diagonal\/} for a Banach algebra $\A$ is an element ${\mathbf M} \in (\A \Tensor \A)^{\ast\ast}$ such that
\[
  a \cdot {\mathbf M} = {\mathbf M} \cdot a \quad\text{and}\qquad a \Delta^{\ast\ast}{\mathbf M} = a \qquad (a \in \A).
\]
\end{definition}
\par
In \cite{Joh2}, Johnson showed that a Banach algebra $\A$ is amenable if and only if it has a virtual diagonal. This allows to introduce a quantified notion of amenability:
\begin{definition}
A Banach algebra $\A$ is called {\it $C$-amenable\/} for some $C \geq 1$ if it has a virtual diagonal of norm at most $C$. The infimum over all $C \geq 1$ such that $\A$ is $C$-amenable is called the
{\it amenability constant\/} of $\A$ and denoted by $\AM_\A$.
\end{definition}
\begin{remark}
It follows from the Alaoglu--Bourbaki theorem (\cite[Theorem V.4.2]{DS}), that the infimum in the definition of $\AM_\A$ is attained, i.e.\ is a minimum.
\end{remark}
\par
Definition \ref{diagdef} has a variant that is better suited for dual Banach algebras. Let $\A$ be a dual Banach algebra with predual $\A_\ast$, and let ${\cal B}_\sigma^2(\A,\comps)$ denote the bounded, bilinear functionals on $\A \times \A$
which are separately $w^\ast$-continuous. Since $\Delta^\ast$ maps $\A_\ast$ into ${\cal B}_\sigma^2(\A,\comps)$, it follows that $\Delta^{\ast\ast}$ drops to an $\A$-bimodule homomorphism $\Delta_\sigma \!: {\cal B}_\sigma^2(\A,\comps)^\ast
\to \A$. We define:
\begin{definition} \label{normdiag}
A {\it normal, virtual diagonal\/} for a dual Banach algebra $\A$ is an element ${\mathrm M} \in {\cal B}_\sigma^2(\A,\comps)^\ast$ such that
\[
  a \cdot {\mathrm M} = {\mathrm M} \cdot a \quad\text{and}\qquad a \Delta_\sigma{\mathrm M} = a \qquad (a \in \A).
\]
\end{definition} 
\begin{remarks}
\item Every dual Banach algebra with a normal, virtual diagonal is Connes-amenable (\cite{CG}).
\item A von Neumann algebra is Connes-amenable if and only if it has a normal, virtual diagonal (\cite{Eff}).
\item The same is true for the measure algebras of locally compact groups (\cite{Run3}).
\end{remarks}
\par
In \cite{Run}, we introduced a stronger variant of Definition \ref{condef} --- called ``strong Connes-amenability'' --- and showed that the existence of a normal, virtual diagonal for a dual Banach algebra was equivalent to it
being strongly Connes-amenable (\cite[Theorem 4.7]{Run}). The following proposition, observed by the late B.\ E.\ Johnson, shows that strong Connes-amenability is even stronger than it seems:
\begin{proposition}
The following are equivalent for a dual Banach algebra $\A$:
\begin{items}
\item There is a normal, virtual diagonal for $\A$.
\item $\A$ has an identity, and every $w^\ast$-continuous derivation from $\A$ into a dual, unital Banach $\A$-bimodule is inner.
\end{items}
\end{proposition}
\begin{proof}
In view of \cite[Theorem 4.7]{Run}, only (i) $\Longrightarrow$ (ii) needs proof. 
\par
Let $E$ be a dual, unital Banach $\A$-bimodule. Due to \cite[Theorem 4.7]{Run}, it is sufficient to show that $D\A \subset E_\sigma$. This, however, is automatically true because
\[
  a \cdot Db = D(ab) - (Da) \cdot b \quad\text{and}\quad (Db) \cdot a = D(ab) - b \cdot Da \qquad (a,b \in \A)
\] 
holds.
\end{proof}
\subsection{Injectivity for Banach modules}
Let $\A$ be a Banach algebra, and let $E$ be a Banach space. Then ${\cal B}(\A,E)$ becomes a left Banach $\A$-bimodule by letting
\[
  (a \cdot T)(x) := T(xa) \qquad (a, x \in \A).
\]
If $E$ is also a left Banach $\A$-module, there is a canonical module homomorphism $\iota \!: E \to {\cal B}(\A,E)$, namely
\[
  \iota(x)a := a \cdot x \qquad (x \in E, \, a \in \A).
\]
\par
For the definition of injective, left Banach modules denote, for any Banach algebra $\A$, by $\A^\#$ the unconditional unitization, i.e.\ we adjoin an identity to $\A$ no matter if $\A$ already has one or not.
Clearly, if $E$ is a left Banach $\A$-module, the module operation extends canonically to $\A^\#$.
\begin{definition} \label{injdef}
Let $\A$ be a Banach algebra. A left Banach $\A$-module $E$ is called {\it injective\/} if $\iota \!: E \to {\cal B}(\A^\#,E)$ has a bounded left inverse which is also a left $\A$-module homomorphism.
\end{definition}
\par
There are various equivalent conditions characterizing injectivity (see, e.g., \cite[Proposition 5.3.5]{LoA}). The following is \cite[Proposition 1.7]{DP}:
\begin{lemma} \label{garthmaksim}
Let $\A$ be a Banach algebra, and let $E$ be a faithful left Banach $\A$-module, i.e.\ if $x \in E$ is such that $a \cdot x =0$ for all $a \in \A$, then $x = 0$. Then $E$ is injective if and only if
$\iota \!: E \to {\cal B}(\A,E)$ has a bounded left inverse which is also an $\A$-module homomorphism.
\end{lemma}
\par
Definition \ref{injdef} and Lemma \ref{garthmaksim} can be adapted to the context of right modules and bimodules in a straightforward way.
\par
The relevance of injectivity in the context of amenable Banach algebras becomes apparent from \cite[Theorem VII.2.20]{Hel} and the duality between injectivity and flatness (\cite[Theorem VII.1.14]{Hel}): A Banach algebra $\A$ with bounded 
approximate identity is amenable if and only if the Banach $\A$-bimodule $\A^\ast$ is injective.
\section{Injectivity of the predual bimodule}
In view of the characterization of amenable Banach algebras just mentioned, one might ask if an analogous statement holds for Connes-amenable, dual Banach algebras $\A$ with $\A^\ast$ replaced by $\A_\ast$. For
von Neumann algebras, this is known to be true (\cite{HelSb}). 
\par
Our first result is true for {\it all\/} dual Banach algebras:
\begin{proposition} \label{DtoC}
Let $\A$ be a dual Banach algebra with identity such that its predual bimodule $\A_\ast$ is injective. Then $\A$ has a normal, virtual diagonal. 
\end{proposition}
\begin{proof}
Consider the short exact sequence
\begin{equation} \label{sequence}
  \{ 0 \} \to \A_\ast \stackrel{\Delta^\ast|_{\A_\ast}}{\to} {\cal B}_\sigma^2(\A,\comps) \to {\cal B}_\sigma^2(\A,\comps)/ \Delta^\ast \A_\ast \to \{ 0 \}.
\end{equation}
Define $P \!: {\cal B}_\sigma^2(\A,\comps) \to \A_\ast$ by letting
\[
  (P\Phi)(a) := \Phi(a,e_\A) \qquad (\Phi \in {\cal B}_\sigma^2(\A,\comps), \, a \in \A),
\]
where $e_\A$ denotes the identity of $\A$. Then it is routinely checked that $P$ is a bounded projection onto $\Delta^\ast \A_\ast$ and thus a left inverse of $\Delta^\ast |_{\A_\ast}$. Hence, (\ref{sequence}) is 
admissible (\cite[Definition 2.3.12]{LoA}).
Since $\A_\ast$ is an injective $\A$-bimodule, there is a bounded $\A$-bimodule homomorphism $\rho \!: {\cal B}_\sigma^2(\A,\comps) \to \A_\ast$ which is a left inverse of $\Delta^\ast |_{\A_\ast}$
(\cite[Proposition 5.3.5]{LoA}). It is routinely checked that $\rho^\ast(e_\A)$ is a normal, virtual diagonal for $\A$.
\end{proof}
\par
As we shall soon see, the converse of Proposition \ref{DtoC} is, in general, false. Nevertheless, for certain $\A$, the injectivity of $\A_\ast$ is indeed equivalent to the existence of a normal virtual diagonal for $\A$
(and even to its Connes-amenability).
\par
We first require a lemma:
\begin{lemma} \label{injlem}
Let $\A$ be a Banach algebra with identity, let $I$ be a closed ideal of $\A$, and let $E$ be a unital Banach $\A$-bimodule such that:
\begin{alphitems}
\item $E$ is injective as a Banach $I$-bimodule.
\item $E$ is faithful both as a left and a right Banach $I$-module.
\end{alphitems} 
Then $E$ is injective as a Banach $\A$-bimodule.
\end{lemma}
\begin{proof}
Turn ${\cal B}(\A \Tensor \A,E)$ into a Banach $\A$-bimodule, by letting
\[
  (a \cdot T)(x \tensor y) := T(x \tensor ya) \quad\text{and}\quad (T \cdot a)(x \tensor y) := T(ax \tensor y) \qquad(a,x,y \in \A).
\]
Define $\iota \!: E \to {\cal B}(\A \Tensor \A,E)$ by letting
\[
  \iota(x)(a \tensor b) := a \cdot x \cdot b \qquad (x \in E, \, a,b \in \A).
\]
Since $\A$ has an identity and $E$ is unital, it is sufficient by (the bimodule analogue of) Lemma \ref{garthmaksim} to show that $\iota$ has a bounded left inverse which is an $\A$-bimodule homomorphism.
\par
By (a), $\iota$ has a bounded left inverse $\rho$ which is an $I$-bimodule homomorphism. We claim that $\rho$ is already an $\A$-bimodule homomorphism. To see this, let $a \in \A$, $T \in {\cal B}(\A \Tensor \A,E)$, and
$b \in I$. Since $I$ is an ideal of $\A$, we obtain that
\[
  b \cdot \rho(a \cdot T) = \rho(ba \cdot T) = ba \cdot \rho(T),
\] 
so that $b \cdot (\rho(a \cdot T) - a \cdot \rho(T)) = 0$. Since $b \in I$ was arbitrary, and since $E$ is a faithful left Banach $I$-module by (b), we obtain $\rho(a \cdot T) = a \cdot \rho(T)$; since $a \in \A$ and
$T  \in {\cal B}(\A \Tensor \A,E)$ were arbitrary, $\rho$ is therefore a left $\A$-module homomorphism. 
\par
Analogously, one shows that $\rho$ is a right $\A$-module homomorphism.
\end{proof}
\par
Our first theorem, considerably improves \cite[Theorem 4.4]{Run}:
\begin{theorem} \label{injthm}
Let $\A$ be an Arens regular Banach algebra which is an ideal in $\A^{\ast\ast}$. Then the following are equivalent:
\begin{items}
\item $\A$ is amenable.
\item $\A^\ast$ is an injective Banach $\A^{\ast\ast}$-bimodule.
\item $\A^{\ast\ast}$ has a normal, virtual diagonal.
\item $\A^{\ast\ast}$ is Connes-amenable.
\end{items}
\end{theorem}
\begin{proof}
(i) $\Longrightarrow$ (ii): We wish to apply Lemma \ref{injlem}. Since $\A$ is amenable, it has a bounded approximate identity. The Arens regularity of $\A$ yields that $\A^{\ast\ast}$ has an identity and that
$\A^\ast$ is a unital Banach $\A^{\ast\ast}$-bimodule. Since $\A$ is amenable and thus a flat $\A$-bimodule over itself (\cite[Theorem VII.2.20]{Hel}), $\A^\ast$ is an injective Banach $\A$-bimodule by (the bimodule version of)
\cite[Theorem VII.1.14]{Hel}. Thus, Lemma \ref{injlem}(a) is satisfied. To see that Lemma \ref{injlem}(b) holds as well, let $\phi \in \A^\ast \setminus \{ 0 \}$. Choose $a \in \A$ such that $\langle a, \phi \rangle \neq 0$.
Let $( e_\alpha )_\alpha$ be a bounded approximate identity for $\A$. Since $\lim_\alpha \langle a e_\alpha, \phi \rangle = \langle a, \phi \rangle \neq 0$, there is $b \in \A$ such that $\langle ab, \phi \rangle =
\langle a, b \cdot \phi \rangle \neq 0$ and thus $b \cdot \phi \neq 0$. Consequently, $\A^\ast$ is faithful as a left Banach $\A$-module. Analogously, one verifies the faithfulness of $\A^\ast$ as a right
Banach $\A$-module.
\par
(ii) $\Longrightarrow$ (iii) is clear by Proposition \ref{DtoC}.
\par
(iii) $\Longrightarrow$ (iv) holds by \cite{CG}.
\par
(iv) $\Longrightarrow$ (i): This is one direction of \cite[Theorem 4.4]{Run}.
\end{proof}
\begin{example}
Let $E$ be a reflexive Banach space with the approximation property, and let $\A$ be ${\cal K}(E)$, the algebra of all compact operators on $E$. Then $\A^\ast$ can be canonically identified with ${\cal N}(E^\ast)$, the
nuclear operators on $E^\ast$, and we have $\A^{\ast\ast} = {\cal B}(E)$. By Theorem \ref{injthm}, we have the equivalence of the following properties:
\begin{items}
\item ${\cal K}(E)$ is amenable.
\item ${\cal N}(E^\ast)$, the space of nuclear operators on $E^\ast$, is an injective Banach ${\cal B}(E)$-bimodule.
\item ${\cal B}(E)$ has a normal, virtual diagonal.
\item ${\cal B}(E)$ is Connes-amenable.
\end{items}
\end{example}
\par
In view of the situation for von Neumann algebras, one might be tempted by Theorem \ref{injthm} to jump to the conclusion that, for a dual Banach algebra $\A$ with predual $\A_\ast$, the injectivity of $\A_\ast$
is equivalent to $\A$ being Connes-amenable or having a normal, virtual diagonal. 
\par
Our next theorem reveals that this is not the case: this gives a negative answer to a question posed by A.\ Ya.\ Helemski\u{\i} 
(\cite[Problem 24]{LoA}). 
\begin{lemma} \label{c0lem}
Let $G$ be a locally compact group, and suppose that ${\cal C}_0(G)$ is injective as a left Banach $M(G)$-module. Then ${\cal C}_0(G)$ is also injective as a left Banach $L^1(G)$-module.
\end{lemma}
\begin{proof}
For $\B = M(G)$ or $\B =L^1(G)$, turn ${\cal B}(\B,{\cal C}_0(G))$ into a left Banach $M(G)$-module in the canonical way, and let $\iota_\B \!: {\cal C}_0(G) \to {\cal B}(\B,{\cal C}_0(G))$ be the respective canonical 
left $M(G)$-module homomorphism.
\par
Since ${\cal C}_0(G)$ is a unital left $M(G)$-module, it is immediate that the homomorphism $\iota_{M(G)} \!: {\cal C}_0(G) \to {\cal B}(M(G),{\cal C}_0(G))$ of left modules, has a linear, bounded left inverse. 
Since ${\cal C}_0(G)$ is injective as a left Banach $M(G)$-module, $\iota_{M(G)}$ has a left inverse $\rho$ which is a bounded homomorphism of left $M(G)$-modules. Let $T \in {\cal B}(M(G),{\cal C}_0(G))$ be such that 
$T |_{L^1(G)} =0$. Since $L^1(G)$ is an ideal in $M(G)$, it follows from the definition of the module action on ${\cal B}(M(G),{\cal C}_0(G))$, that $f \cdot T = 0$ for all $f \in L^1(G)$ and therefore 
$f \cdot \rho(T) = \rho(f \cdot T) = 0$ for all $f \in L^1(G)$. Since ${\cal C}_0(G)$ is a faithful left $L^1(G)$-module, this means that $\rho(T) = 0$. Since $L^1(G)$ is complemented in $M(G)$, it follows 
that $\rho \!: {\cal B}(M(G),{\cal C}_0(G)) \to {\cal C}_0(G)$ drops to bounded homomorphism of left $M(G)$-modules $\tilde{\rho} \!: {\cal B}(L^1(G),{\cal C}_0(G)) \to {\cal C}_0(G)$, which is easily seen to be a left inverse of 
$\iota_{L^1(G)}$.
\par
Since $\tilde{\rho}$ is trivially a homomorphism of left $L^1(G)$-modules, Lemma \ref{garthmaksim} yields the injectivity of ${\cal C}_0(G)$ as a left Banach $L^1(G)$-module.
\end{proof}
\begin{theorem} \label{BnotC}
Let $G$ be a locally compact group. Then ${\cal C}_0(G)$ is an injective Banach $M(G)$-bimodule if and only if $G$ is finite.
\end{theorem}
\begin{proof}
Suppose that ${\cal C}_0(G)$ is an injective Banach $M(G)$-bimodule. By \cite[Proposition VII.2.1]{Hel} and Lemma \ref{c0lem}, ${\cal C}_0(G)$ is also injective as a left Banach $L^1(G)$-module. By \cite[Theorem 3.8]{DP},
this means that $G$ must be finite.
\par
The converse is obvious.
\end{proof}
\begin{remark}
In contrast, it was proven in \cite{Run2}, for a locally compact group $G$, that $M(G)$ is Connes-amenable --- and, equivalently, has a normal, virtual diagonal by \cite{Run3} --- if and only if $G$ is 
amenable.
\end{remark}
\section{Fourier--Stieltjes algebras of locally compact groups}
The Fourier--Stieltjes algebra $B(G)$ of a locally compact group $G$ was introduced by P.\ Eymard in \cite{Eym} along with the Fourier algebra $A(G)$. We refer to \cite{Eym} for further information on these algebras. It is straightforward to 
see that $B(G)$ is a dual Banach algebra --- with predual $\cstar(G)$ --- for any locally compact group $G$ whereas $A(G)$ need not even be a dual space (unless $G$ is compact, of course).
\par
Let $G$ be a locally compact group $G$ with an abelian subgroup of finite index. Then $A(G)$ is amenable and $w^\ast$-dense in $B(G)$, so that $B(G)$ is Connes-amenable. In fact, a formally stronger conclusion holds:
\begin{proposition} \label{FSTprop}
Let $G$ be a locally compact group with an abelian subgroup of finite index. Then $B(G)$ has a normal, virtual diagonal.
\end{proposition}
\begin{proof}
Let $H$ be a an abelian subgroup of $G$ such that $n := [G:H] < \infty$. Replacing $H$ by its closure, we may suppose that $H$ is closed and thus open. Consequently, the restriction map
from $B(G)$ onto $B(H)$ is surjective so that
\[
  B(G) \cong B(H)^n \cong M(\hat{H})^n,
\]
where $\hat{H}$ is the dual group of $H$. By \cite{Run3}, $M(\hat{H})$ has a normal, virtual diagonal. It is easy to see that therefore $M(\hat{H})^n \cong B(G)$ must have a normal, virtual diagonal
as well.
\end{proof}
\par
In view of \cite[Theorem 5.2]{Run4}, we conjecture that the converse of Proposition \ref{FSTprop} holds as well --- even with the existence of a normal, virtual diagonal replaced by Connes-amenability. We have, however, been 
unable to confirm this conjecture for arbitrary locally compact groups. In the remainder of this section, we shall prove partial converses of Proposition \ref{FSTprop} for groups with certain additional properties.
\par
Given a family $( G_\alpha )_\alpha$ of locally compact groups, we denote by $\prod_\alpha G_\alpha$ its direct product equipped with the product topology. 
\begin{lemma} \label{FSTlem1}
Let $( G_\alpha )_\alpha$ be a family of locally compact groups, let $G := \prod_\alpha G_\alpha$, and let $\pi_\alpha \!: B(G) \to B(G_\alpha)$ be the canonical quotient map for each index $\alpha$. Then
\[
  \pi := \bigoplus_\alpha\pi_\alpha \!: B(G) \to \text{$\ell^\infty$-}\bigoplus_\alpha B(G_\alpha).
\] 
is a $w^\ast$-continuous algebra homomorphism with $w^\ast$-dense range.
\end{lemma}
\begin{proof}
Since each $G_\alpha$ --- viewed as a subgroup of $G$ --- is open, it follows that each map $\pi_\alpha$ is $w^\ast$-continuous. Consequently, $\pi$ is $w^\ast$-continuous.
\par
We may view each $B(G_\alpha)$ as a closed subalgebra of $\text{$\ell^\infty$-}\bigoplus_\alpha B(G_\alpha)$ in a canonical fashion. To establish that $\pi$ has $w^\ast$-dense range, it is sufficent to show that $B(G_\alpha) \subset\pi(B(G))$
for each index $\alpha$. Fix $\alpha$, and let $\chi_\alpha \!: G \to \comps$ denote the indicator funtion of $G_\alpha$. Since $G_\alpha$ is an open subgroup of $G$, we have that $\chi_\alpha \in B(G)$. Clearly, $\pi$ maps
$\chi_\alpha B(G)$ onto $B(G_\alpha)$.
\end{proof}
\begin{lemma} \label{FSTlem2}
Let $( G_\alpha )_\alpha$ be a family of finite groups, let $G := \prod_\alpha G_\alpha$, and suppose that $B(G)$ is Connes-amenable. Then $\sup_\alpha \AM_{B(G_\alpha)}$ is finite.
\end{lemma}
\begin{proof}
Since $B(G)$ is Connes-amenable, the same is true for $\text{$\ell^\infty$-}\bigoplus_\alpha B(G_\alpha)$ by \cite[Proposition 4.2(ii)]{Run}. Since each group $G_\alpha$ is finite, $B(G_\alpha)$ is finite-dimensional so that
\[
  \text{$\ell^\infty$-}\bigoplus_\alpha B(G_\alpha) = \left( \text{$c_0$-}\bigoplus_\alpha B(G_\alpha) \right)^{\ast\ast}.
\]
Since $\A := \text{$c_0$-}\bigoplus_\alpha B(G_\alpha)$ is an ideal in $\text{$\ell^\infty$-}\bigoplus_\alpha B(G_\alpha)$, it follows from \cite[Theorem 4.4]{Run} that $\A$ is amenable. Since amenability constants only shrink under
passage to quotients, we conclude that $\AM_{B(G_\alpha)} \leq \AM_\A$ holds for each index $\alpha$. 
\end{proof}
\par
We can now prove our first partial converse of Proposition \ref{FSTprop}:
\begin{theorem} \label{compactthm}
Let $( G_\alpha )_\alpha$ be a family of finite groups, and let $G := \prod_\alpha G_\alpha$. Then the following are equivalent for $G$:
\begin{items}
\item All but finitely many of the groups $G_\alpha$ are abelian.
\item $A(G)$ is amenable.
\item $B(G)$ has a normal, virtual diagonal.
\item $B(G)$ is Connes-amenable.
\end{items}
\end{theorem}
\begin{proof}
(i) $\Longrightarrow$ (ii) is well known (\cite[Theorem 4.5]{JohF} or \cite[Corollary 4.3]{LLW}).
\par
(ii) $\Longrightarrow$ (iii): Since $G$ is compact, we have $B(G) = A(G)$, so that $B(G)$ is amenable and thus has a virtual diagonal $\mathbf{M} \in (B(G) \Tensor B(G))^{\ast\ast}$. Restricting $\mathbf{M}$ to
${\cal B}_\sigma^2(B(G),\comps) \subset (B(G) \Tensor B(G))^\ast$, we obtain a normal, virtual diagonal for $B(G)$.
\par
(iii) $\Longrightarrow$ (iv) is clear.
\par
(iv) $\Longrightarrow$ (i): Assume that there is a subfamily $( G_{\alpha_n} )_{n=1}^\infty$ of $( G_\alpha )_\alpha$, such that $G_{\alpha_n}$ is not abelian for each $n \in \posints$. For $n \in \posints$, define
\[
  H_n := G_{\alpha_\frac{n(n+1)}{2}} \times \cdots \times G_{\alpha_{\frac{n(n+1)}{2}+(n-1)}}.
\]
Let $H := \prod_{n=1}^\infty H_n$. Since the restriction map from $B(G)$ to $B(H)$ is a $w^\ast$-continuous algebra homomorphism with $w^\ast$-dense range (even surjective), \cite[Proposition 4.2(ii)]{Run} shows that $B(H)$ is
Connes-amenable as well. It therefore follows from Lemma \ref{FSTlem1} that $\sup_{n \in \posints} \AM_{B(H_n)} < \infty$. This, however, contradicts \cite[Corollary 4.2 and Proposition 4.3]{JohF} which assert that
$\AM_{B(H_n)} \geq \left( \frac{3}{2} \right)^n$ for each $n \in \posints$. 
\end{proof}
\par
The groups considered in Theorem \ref{compactthm} are compact. For amenable, discrete groups, another partial converse of Proposition \ref{FSTprop} holds:
\begin{theorem} \label{discrthm}
The following are equivalent for an amenable, discrete group $G$:
\begin{items}
\item $G$ has an abelian subgroup of finite index.
\item $A(G)$ is amenable.
\item $B(G)$ has a normal, virtual diagonal.
\item $B(G)$ is Connes-amenable.
\end{items}
\end{theorem}
\begin{proof}
(i) $\Longrightarrow$ (iii) is Proposition \ref{FSTprop} and (iii) $\Longrightarrow$ (iv) is clear.
\par
(iv) $\Longrightarrow$ (ii): Let $E$ be a Banach $A(G)$-bimodule, and let $D \!: A(G) \to E^\ast$ be a bounded derivation. Since $G$ is amenable, $A(G)$ has a bounded approximate identity by Leptin's
theorem (\cite[Theorem 7.1.3]{LoA}). Hence, by \cite[Proposition 2.1.5]{LoA}, we may suppose that $E$ is pseudo-unital (\cite[Definition 2.1.4]{LoA}). Let $\tau$ denote the multiplier topology on $B(G)$, i.e.\ 
a net $( f_\alpha )_\alpha$ in $B(G)$ converges to $f \in B(G)$ with respect to $\tau$ if and only if $\| f_\alpha g - f g \|_{A(G)} \to 0$ for each $g \in A(G)$. By \cite[Proposition 2.1.6]{LoA}, the module actions of
$A(G)$ on $E$ extend to $B(G)$ in a canonical manner; it is immediate that, for $x \in E$, the module actions
\[
  B(G) \to E, \quad f \mapsto \left\{ \begin{array}{c} f \cdot x, \\ x \cdot f \end{array} \right.
\]
are $\tau$-norm-continuous. Furthermore, \cite[Proposition 2.1.6]{LoA} asserts that $D$ extends to a $\tau$-$w^\ast$-continuous derivation $\tilde{D} \!: B(G) \to E^\ast$. Since $G$ is discrete and since $A(G)$ is
regular, $\tau$ and the $w^\ast$-topology of $B(G)$ coincide on norm bounded subsets. From the Krein--\v{S}mulian theorem (\cite[Theorem V.5.7]{DS}), we conclude that $E^\ast$ is a normal, dual Banach $B(G)$-module and
that $\tilde{D}$ is $w^\ast$-continuous. Consequently, $\tilde{D}$ is inner, and so is $D$.
\par
(ii) $\Longrightarrow$ (i) is the difficult direction of \cite[Theorem 5.2]{Run4}.
\end{proof}
\par
The {\it reduced Fourier--Stieltjes algebra\/} $B_r(G)$ of a locally compact group $G$, was also introduced in \cite{Eym}. It is the dual of the reduced group $\cstar$-algebra $C^\ast_r(G)$ and is a $^\ast$-closed ideal in $B(G)$.
As another consequence of Theorem \ref{discrthm}, we obtain (compare \cite[Theorem 4.4]{RS}):
\begin{corollary}
Let $G$ be a discrete group. Then $B_r(G)$ is Connes-amenable if and only if $G$ has an abelian subgroup of finite index.
\end{corollary}
\begin{proof}
Suppose that $B_r(G)$ is Connes-amenable. Then $B_r(G)$ has an identity (\cite[Proposition 4.1]{Run}) and thus equals $B(G)$. Consequently, $G$ is amenable.
\par
The rest is a straightforward consequence of Theorem \ref{discrthm}.
\end{proof}
\section{Weak almost periodicity, $w^\ast$-weak continuity, and normality}
One of the unsatisfactory sides of dealing with Connes-amenability for dual Banach algebras is the apparent lack of a suitable intrinsic characterization in terms of virtual diagonals. Dual Banach algebras with a normal, virtual
diagonal are Connes-amenable, but the converse is likely to be false in general. 
For von Neumann algebras (\cite{Eff}) and measure algebras (\cite{Run3}), (A) and (B) are equivalent, but in both cases the methods employed to prove this
equivalence give no clue about how to tackle the general case.
\par
In this section, we pursue a different approach towards a ``virtual diagonal characterization'' for Connes-amenable, dual Banach algebras. The main problem with trying to prove that (A) implies (B) is that, for a general dual Banach algebra $\A$,
the module ${\cal B}_\sigma^2(\A,\comps)^\ast$ need not be normal. In this section, we show that every dual Banach $\A$-bimodule has what one might call a largest normal quotient. Using this idea to modify the definition of a normal, virtual
diagonal, we then obtain the desired characterization (Theorem \ref{diagthm} below).
\par
We begin with recalling the notion of weak almost periodicity (in a slightly more general context than usual):
\begin{definition} \label{wapdef}
Let $\A$ be a Banach algebra, and let $E$ be a Banach $\A$-bimodule. Then an element $x \in E$ is called {\it weakly almost periodic\/} if the module maps
\[
  \A \to E, \quad a \mapsto  \left\{ \begin{array}{c} a \cdot x, \\ x \cdot a \end{array} \right.
\]
are weakly compact. The collection of all weakly almost periodic elements of $E$ is denoted by $\WAP(E)$.
\end{definition}
\begin{remark}
It follows easily from Grothendieck's double limit criterion that
\[
  \WAP(\A^\ast) = \{ \phi \in \A^\ast : \text{$\A \ni a \mapsto a \cdot \phi$ is weakly compact} \},
\]
so that our choice of terminology is consistent with the usual one as used in \cite{LL}, for instance.
\end{remark}
\par
The reason why we are interested in weak almost periodicity in the context of dual Banach algebras is that it is closely related to the normality of dual Banach modules:
\begin{proposition} \label{wap1}
Let $\A$ be a dual Banach algebra, and let $E$ be a Banach $\A$-bimodule such that $E^\ast$ is normal. Then $E= \WAP(E)$ holds.
\end{proposition}
\begin{proof}
Let $x \in E$. Since $E^\ast$ is normal, it follows immediately from the definition of $\sigma(E,E^\ast)$ that the maps
\begin{equation} \label{modops}
  \A \to E, \quad a \mapsto  \left\{ \begin{array}{c} a \cdot x, \\ x \cdot a \end{array} \right.
\end{equation}
are $\sigma(\A,\A_\ast)$-$\sigma(E,E^\ast)$-continuous, where $\A_\ast$ is the predual bimodule of $\A$. Since the closed unit ball of $\A$ is $\sigma(\A,\A_\ast)$-compact by the
Alaoglu--Bourbaki theorem, it follows that the maps (\ref{modops}) are weakly compact.
\end{proof}
\par
At first glance, one might conjecture that the converse of Proposition \ref{wap1} holds as well. This, however, is not true:
\begin{example}
Let $G$ be a locally compact group. Recall that a continuous, bounded function on $G$ is called weakly almost periodic if its orbit under left (and, equivalently, under right) translation is weakly compact (\cite{Bur}).
We denote the space of all weakly almost periodic functions on $G$ by $\WAP(G)$. It is easy to see that $\WAP(G)$ is a commutative $\cstar$-algebra. Its character space carries a natural semigroup structure (with separately
continuous multiplication) that extends multiplication on $G$ (see \cite{Bur} for details): we denote this so-called weakly almost periodic compactification of $G$ by $wG$. Identifying $\WAP(G)^\ast$ with $M(wG)$, we can equip
$\WAP(G)^\ast$ with a convolution type product turning it into a dual Banach algebra. Via integration, the dual Banach algebra $M(G)$ can be identified with a subalgebra of $\WAP(G)^\ast$, so that $E = \WAP(G)$ becomes a
a Banach $M(G)$-bimodule in a canonical fashion. It is straightforward to check that $\WAP(E) = \WAP(G)$ whereas $E^\ast$ is normal as a dual Banach $M(G)$-bimodule if and only if $G$ is compact.
\end{example}
\par
We therefore feel justified to introduce a new definition:
\begin{definition}
Let $\A$ be a dual Banach algebra with predual $\A_\ast$, and let $E$ be a Banach $\A$-bimodule. Then an element $x \in E$ is called {\it $w^\ast$-weakly continuous\/} if the module maps
\[
  \A \to E, \quad a \mapsto  \left\{ \begin{array}{c} a \cdot x, \\ x \cdot a \end{array} \right.
\]
are $\sigma(\A,\A_\ast)$-$\sigma(E,E^\ast)$-continuous. The collection of all $w^\ast$-weakly continuous elements of $E$ is denoted by $\wWC(E)$.
\end{definition}
\begin{remarks}
\item It is easy to see that $\wWC(E)$ is a closed submodule of $E$.
\item If $F$ is another Banach $\A$-bimodule and if $\theta \!: E \to F$ is a bounded $\A$-bimodule homomorphism, then $\theta(\wWC(E)) \subset \wWC(F)$ holds.
\item It is implicit in the proof of Proposition \ref{wap1} that $\wWC(E) \subset \WAP(E)$.
\item For a locally compact group $G$, $\A = M(G)$, and $E = \WAP(G)$, we have that $\WAP(E) = \WAP(G)$ whereas $\wWC(E) = {\cal C}_0(G)$. Hence, $\wWC(E) \subsetneq \WAP(E)$ holds whenever $G$ is not compact. 
\end{remarks}
\begin{proposition} \label{wap2}
Let $\A$ be a dual Banach algebra, and let $E$ be a Banach $\A$-bimodule. Then then following are equivalent:
\begin{items}
\item $E^\ast$ is normal.
\item $E = \wWC(E)$.
\end{items}
\end{proposition}
\begin{proof}
(i) $\Longrightarrow$ (ii): The argument used to prove Proposition \ref{wap1} does in fact yield the stronger assertion (ii).
\par
(ii) $\Longrightarrow$ (i): This is proved in the same way as (i) $\Longrightarrow$ (ii), only with the r\^oles of $E$ and $E^\ast$ interchanged.
\end{proof}
\par
It follows from Proposition \ref{wap2} that, for any Banach $\A$-bimodule $E$, the dual module $\wWC(E)^\ast$ is normal. We therefore obtain:
\begin{corollary}
A dual Banach algebra $\A$ is Connes-amenable if, for every Banach $\A$-bimodule $E$, every $w^\ast$-continuous derivation $D \!: \A \to \wWC(E)^\ast$ is inner.
\end{corollary}
\par
Since any dual Banach algebra is a normal dual Banach module over itself, we obtain as another consequence of Proposition \ref{wap2}:
\begin{corollary} \label{dualcor}
Let $\A$ be a dual Banach algebra with predual bimodule $\A_\ast$. Then $\A_\ast \subset \wWC(\A^\ast)$ holds.
\end{corollary}
\par
Let $\A$ be a dual Banach algebra with predual $\A_\ast$, and let $\Delta \!: \A \Tensor \A \to \A$ be the multiplication map. From Corollary \ref{dualcor}, we conclude that $\Delta^\ast$ maps
$\A_\ast$ into $\wWC((\A \Tensor \A)^\ast)$. Consequently, $\Delta^{\ast\ast}$ drops to homomorphism $\Delta_{\wWC} \!: \wWC((\A \Tensor \A)^\ast)^\ast \to \A$. 
\par
With these preparations made, we can now characterize the Connes-amenable, dual Banach algebras through the existence of certain virtual diagonals:
\begin{definition} \label{wapdiag}
Let $\A$ be a dual Banach algebra. A {\it $\wWC$-virtual diagonal\/} for $\A$ is an element $\mathrm{M} \in \wWC((\A \Tensor \A)^\ast)^\ast$ such that
\[
  a \cdot \mathrm{M} = \mathrm{M} \cdot a \quad\text{and}\quad a \Delta_{\wWC} \mathrm{M} =a \qquad (a \in \A).
\]
\end{definition}
\begin{theorem} \label{diagthm}
The following are equivalent for a dual Banach algebra $\A$:
\begin{items}
\item $\A$ is Connes-amenable.
\item There is a $\wWC$-virtual diagonal for $\A$.
\end{items}
\end{theorem}
\begin{proof}
(i) $\Longrightarrow$ (ii): First, note that $\A \Tensor \A$ is canonically mapped into $\wWC((\A \Tensor \A)^\ast)^\ast$; in order to make notation not more complicated than necessary (and than it already is),
we just write those elements of $\wWC((\A \Tensor \A)^\ast)^\ast$ that lie in the canonical image of $\A \Tensor \A$ as tensors. 
\par
By \cite[Proposition 4.1]{Run}, $\A$ has an identity $e_\A$. Define a derivation
\[
  D \!: \A \to \wWC((\A \Tensor \A)^\ast)^\ast, \quad a \mapsto a \tensor e_\A - e_\A \tensor a ( = a \cdot(e_\A \tensor e_\A) - (e_\A \tensor e_\A) \cdot a).
\]
Since the dual module $\wWC((\A \Tensor \A)^\ast)^\ast$ is normal, it follows that $D$ is $w^\ast$-continuous. Clearly, $D$ attains its values in the $w^\ast$-closed submodule $\ker \Delta_{\wWC}$.
Hence, there is $\mathrm{N} \in \ker \Delta_{\wWC}$ such that
\[
  D a = a \cdot \mathrm{N} - \mathrm{N} \cdot a \qquad (a \in \A).
\]
Letting $\mathrm{M} := e_\A \tensor e_\A - \mathrm{N}$, we obtain an element as required by Definition \ref{wapdiag}.
\par
(ii) $\Longrightarrow$ (i): Clearly, (ii) implies that $\A$ has an identity. Let $E$ be a normal, dual Banach $\A$-bimodule --- which we may suppose without loss of generality to be unital ---, and let
$D \!: \A \to E$ be a $w^\ast$-continuous derivation. Define
\[
  \theta_D \!: \A \Tensor \A \to E, \quad a \tensor b \mapsto a \cdot Db.
\]
By Lemma \ref{waplem} below, $\theta_D^\ast$ maps the predual module $E_\ast$ of $E$ into $\wWC((\A \Tensor \A)^\ast)$. Hence, $\Theta_D := (\theta_D^\ast |_{E_\ast} )^\ast$ maps
$\wWC((\A \Tensor \A)^\ast)^\ast$ into $E$. 
\par
Let $\mathrm{M} \in \wWC((\A \Tensor \A)^\ast)^\ast$ be a $\wWC$-virtual diagonal for $\A$, and let $x := \Theta_D(\mathrm{M})$. A more or less verbatim copy of the argument given
in the proof of \cite[Theorem 2.2.4]{LoA} then shows that $D$ is the inner derivation implemented by $x$.
\end{proof}
\par
To complete the proof of Theorem \ref{diagthm}, we require the technical Lemma \ref{waplem} below. To make its proof more transparent, we introduce new notation: Given a dual Banach algebra $\A$ and a left Banach $\A$-module $E$, 
we define
\[
  \wWC_l(E) := \{ x \in E : \text{$\A \ni a \mapsto a \cdot x$ is $w^\ast$-weakly continuous} \};
\]
similarly, we define $\wWC_r(E)$ for a right Banach $\A$-module $E$. Clearly, $\wWC_l(E) \cap \wWC_r(E)$ equals $\wWC(E)$ if $E$ is a Banach $\A$-bimodule.
\begin{lemma} \label{waplem}
Let $\A$ be a dual Banach algebra with identity, let $E$ be a normal, dual Banach $\A$-bimodule with predual bimodule $E_\ast$, and let $D \!: \A \to E$ be a $w^\ast$-continuous derivation. Then the adjoint of
\[
  \theta_D \!: \A \Tensor \A \to E, \quad a \tensor b \mapsto a \cdot Db.
\]
maps $E_\ast$ into $\WAP((\A \Tensor \A)^\ast)^\ast$.
\end{lemma}
\begin{proof}
Clearly, $\theta_D$ is a homomorphism of left Banach $\A$-modules, so that $\theta_D^\ast$ is a homomorphism of right Banach $\A$-modules. It follows from Proposition \ref{wap2} that
\[ 
  \theta_D^\ast(E_\ast) \subset \theta_D^\ast(\wWC(E^\ast)) \subset \theta_D^\ast(\wWC_r(E^\ast)) \subset \wWC_r((\A \Tensor \A)^\ast).
\]
\par
Let $e_\A$ denote the identity of $\A$, and let $R := e_\A \tensor \A$. Then $R$ is a closed submodule of the right Banach $\A$-module $\A \Tensor \A$ such that we have a direct sum $\A \Tensor \A = \ker \Delta \oplus R$ of
right Banach $\A$-modules. Consequently, we have a direct sum $(\A \Tensor \A)^\ast = (\ker \Delta)^\ast \oplus R^\ast$ of left Banach $\A$-modules. 
Let $\theta_1 := \theta_D |_{\ker \Delta}$ and $\theta_2 := \theta_D |_{R}$, so that $\theta_D = \theta_1 \oplus \theta_2$
and, consequently, $\theta_D^\ast = \theta_1^\ast \oplus \theta_2^\ast$. It is easy to check that $\theta_1^\ast$ is a homomorphism of left Banach $\A$-modules, so that $\theta_1^\ast(E_\ast) \subset \wWC_l((\A \Tensor \A)^\ast)$.
Clearly, the right Banach $\A$-modules $R$ and $\A$ are canonically isomorphic. Identifying, $R$ with $\A$, we see that $\theta_2^\ast$ is nothing by the pre-adjoint of the $w^\ast$-continuous linear map $D$. In view of Corollary 
\ref{dualcor}, it follows that
\[
  \theta_2^\ast(E_\ast) \subset \A_\ast \subset \wWC_l(\A^\ast) = \wWC_l(R^\ast) \subset \wWC_l((\A \Tensor \A)^\ast).  
\]
All in all, $\theta_D^\ast(E_\ast) \subset \wWC_l((\A \Tensor \A)^\ast)$ holds.
\end{proof}
\par
Finally, we return to weak almost periodicity in the sense of Definition \ref{wapdef}.
\par
Recall (from \cite{LL}, for instance) the notion of a left introverted subspace of the dual of a Banach algebra: A right $\A$-submodule $E$ of $\A^\ast$ is called {\it left introverted\/} if, for any $\phi \in E$ and $A \in E^\ast$ 
the functional $A \cdot \phi \in \A^\ast$ defined by
\[
  \langle a, A \cdot \phi \rangle := \langle \phi \cdot a,A \rangle \qquad (a \in \A).
\]
lies again in $E$. This can be used to turn $E^\ast$ into a Banach algebra by letting
\begin{equation} \label{intro}
  \langle \phi, AB \rangle := \langle B \cdot \phi, A \rangle \qquad \qquad (A, B \in E^\ast, \, \phi \in E). 
\end{equation}
\par
We can use this construction to define, for an arbitrary Banach algebra, a dual Banach algebra with a certain universality property:
\begin{theorem}
Let $\A$ be a Banach algebra. Then $\WAP(\A^\ast)$ is a left introverted subspace of $\A^\ast$ such that $\WAP(\A^\ast)^\ast$, equipped with the product defined in\/ {\rm (\ref{intro})}, is a dual 
Banach algebra with the following universal property: Whenever $\B$ is a dual Banach algebra, and $\theta \!: \A \to \B$ is a bounded algebra homomorphism, then there is a unique $w^\ast$-continuous algebra homomorphism
$\pi \!: \WAP(\A^\ast)^\ast \to \B$ such that the diagram
\begin{equation} \label{diagram}
  \begin{diagram}[midshaft,middle,tight]
    \A & \rTo^\iota        & \WAP(\A^\ast)^\ast \\
       & \rdTo_\theta      & \dTo_{\pi} \\
       &                   & \B
  \end{diagram}
\end{equation}
commutes where $\iota \!:\A \to \WAP(\A^\ast)^\ast$ is the canonical map.
\end{theorem}
\begin{proof}
By \cite[Lemma 1.2]{LL}, $\WAP(\A^\ast)$ is left introverted, and by \cite[Lemma 1.4]{LL}, $\WAP(\A^\ast)^\ast$ is a dual Banach algebra (the commutativity hypothesis from that lemma is not required for this
particular assertion).
\par
Let $\B$ be a dual Banach algebra with predual $\B_\ast$, and let $\theta \!: \A \to \B$ be a bounded algebra homomorphism. By Corollary \ref{dualcor}, $\B_\ast \subset \wWC(\B^\ast) \subset \WAP(\B^\ast)$ holds.
Furthermore, it is easy to see that $\theta^\ast(\WAP(\B^\ast)) \subset \WAP(\A^\ast)$. Letting $\pi := (\theta^\ast |_{\B_\ast})^\ast$, we obtain a $w^\ast$-continuous map $\pi \!: \WAP(\A^\ast)^\ast \to \B$ such that 
(\ref{diagram}) commutes; since $\iota(\A)$ is $w^\ast$-dense in $\WAP(\A^\ast)^\ast$, this uniquely determines $\pi$. Clearly, $\pi$ is multiplicative if restricted to $\iota(\A)$. Since 
multiplication in both $\WAP(\A^\ast)^\ast$ and $\B$ is separately $w^\ast$-continuous, the $w^\ast$-density of $\iota(\A)$ and the $w^\ast$-continuity of $\pi$ show that $\pi$ is in fact multiplicative on all of $\WAP(\A^\ast)^\ast$.
\end{proof}
\begin{remarks}
\item If $\A$ is Arens regular, $\WAP(\A^\ast) = \A^\ast$ holds, so that $\WAP(\A^\ast)^\ast$ is nothing but $\A^{\ast\ast}$ equipped with either Arens product. 
\item There is no need for $\iota \!: \A \to \WAP(\A^\ast)^\ast$ to be injective, let alone an isometry: If $E$ is a non-reflexive Banach space with the approximation property, we have
$\WAP({\cal K}(E)^\ast) = \{ 0 \}$ by \cite[Theorem 3]{You} (and, consequently, $\WAP({\cal K}(E)^\ast)^\ast = \{ 0 \}$).
\item If $\A$ is a dual Banach algebra, then $\A$ embeds isometrically into $\WAP(\A^\ast)^\ast$ by Corollary \ref{dualcor}.
\end{remarks}
\par
Since the image of $\A$ in $\WAP(\A^\ast)^\ast$ is $w^\ast$-dense for any Banach algebra $\A$, the amenability of $\A$ immediately yields the Connes-amenability of $\WAP(\A^\ast)^\ast$. The converse is clearly false: 
\begin{example}
The dual space ${\cal B}(\ell^2)$ of $E := \ell^2 \Tensor \ell^2$ lacks the approximation property (\cite{Sza}). Consequently, ${\cal K}(E)$ cannot have a bounded approximate identity, let alone be amenable (\cite[Corollary 3.1.5]{LoA}).
Since $E$ is not reflexive, $\WAP({\cal K}(E)^\ast)^\ast = \{ 0 \}$ is trivially Connes-amenable.
\end{example}
\par
Nevertheless, for certain Banach algebras $\A$ the Connes-amenability of $\WAP(\A^\ast)^\ast$ is indeed equivalent to the amenability of $\A$:
\begin{proposition}
The following are equivalent for a locally compact group $G$:
\begin{items}
\item $G$ is amenable.
\item $L^1(G)$ is amenable.
\item $\WAP(L^1(G)^\ast)^\ast$ is Connes-amenable.
\end{items}
\end{proposition}
\begin{proof}
(i) $\Longleftrightarrow$ (ii) is \cite[Theorem 2.5]{Joh1}, and (ii) $\Longrightarrow$ (iii) is clear by \cite[Proposition 4.2(i)]{Run}.
\par
(iii) $\Longrightarrow$ (i): Suppose that $\WAP(L^1(G)^\ast)^\ast$ is Connes-amenable. First, note that the dual Banach algebras $\WAP(L^1(G)^\ast)^\ast$ and $\WAP(G)^\ast$ (mentioned earlier) are identical. 
Since ${\cal C}_0(G) \subset \WAP(G)$, restriction is a $w^\ast$-continuous algebra homomorphism from $\WAP(G)^\ast$ onto $M(G)$, so that $M(G)$ is Connes-amenable by \cite[Propositon 4.2(ii)]{Run}. By \cite{Run2}, this means that 
$G$ is amenable.
\end{proof}
\par
Let $G$ be a compact group. By \cite{Run3}, $\WAP(G)^\ast = M(G)$ then has a normal, virtual diagonal. We do not know if this is still true for locally compact, but non-compact $G$. We suspect, but have been unable to prove, 
that $\WAP(G)^\ast$ with $G$ amenable, but not compact, is an example of a Connes-amenable, dual Banach algebra which fails to have a normal, virtual diagonal.
\dated
\vfill
\begin{tabbing}
{\it Author's address\/}: \= Department of Mathematical and Statistical Sciences \\
\> University of Alberta \\
\> Edmonton, Alberta \\
\> Canada T6G 2G1 \\[\medskipamount]
{\it E-mail\/}: \> {\tt vrunde@ualberta.ca} \\[\medskipamount]
{\it URL\/}: \> {\tt http://www.math.ualberta.ca/$^\sim$runde/} 
\end{tabbing} 
\end{document}